\documentclass[letterpaper, 10pt, conference]{ieeeconf}
\IEEEoverridecommandlockouts                              
\usepackage[dvipsnames]{xcolor}
\usepackage[utf8]{inputenc}
\usepackage{epsfig}
\usepackage{textcomp}
\usepackage{comment}
\usepackage{graphicx}
\usepackage{amssymb}
\usepackage{gensymb}
\usepackage{amsmath}

\usepackage{algorithm}
\usepackage{algpseudocode}
\usepackage[labelfont={bf},font=small]{caption}
\usepackage[none]{hyphenat}
\usepackage{mathtools, cuted}
\usepackage[noadjust, nobreak]{cite}
\usepackage{tabularx}
\usepackage{float}
\usepackage{pifont}
\usepackage[]{placeins}
\usepackage{makecell}
\usepackage{tikz}
\usepackage[framemethod=tikz]{mdframed}
\usepackage{afterpage}
\usepackage{stfloats}
\usepackage{atbegshi}
\usepackage{acronym}
\usepackage{draftwatermark}

\makeatletter
\let\NAT@parse\undefined
\makeatother
\usepackage{hyperref}
\hypersetup{
colorlinks=true, 
breaklinks=true, 
urlcolor= black, 
citecolor=black,
linkcolor= blue, 
bookmarksopen=true,
pdftitle={RNN-based linear parameter varying adaptive model predictive control for autonomous driving}, 
pdfauthor={Y. Kebbati},
}

\title{\LARGE \bf
RNN-based linear parameter varying adaptive model predictive control for autonomous driving
}

\include{acronyms}

\author{Yassine Kebbati$^{1*}$, Naima Ait-Oufroukh$^{1}$, Dalil Ichalal$^{1}$ and Vincent Vigneron$^{1,2}$ 
\thanks{$^*$Corresponding author: \href{mailto:yassine.kebbati@univ-evry.fr}{yassine.kebbati@univ-evry.fr}}
\thanks{$^1$IBISC-EA4526, univ Evry, université Paris-Saclay, France}
\thanks{$^2$ School of Applied Sciences (FCA),UNICAMP, Limeira, Brazil}
}

\makeatletter
\newcommand*{\rom}[1]{\expandafter\@slowromancap\romannumeral #1@}
\makeatother

\begin{document}

\maketitle
\thispagestyle{empty}
\pagestyle{empty}


\SetWatermarkText{Preprint}

\begin{abstract}                          
Autonomous driving is a complex and highly dynamic process that ensures controlling the coupled longitudinal and lateral vehicle dynamics. Model predictive control, distinguished by its predictive feature, optimal performance, and ability to handle constraints, makes it one of the most promising tools for this type of control application. The content of this article handles the problem of autonomous driving by proposing an adaptive linear parameter varying model predictive controller (LPV-MPC), where the controller's prediction model is adaptive by means of a recurrent neural network. The proposed LPV-MPC is further optimized by a hybrid Genetic and Particle Swarm Optimization Algorithm (GA-PSO). The developed controller is tested and evaluated on a challenging track under variable wind disturbance.

\end{abstract}




                               

\section{INTRODUCTION}

The ever-increasing number of vehicles is inducing terrible traffic conditions and increasing air pollution in today's world. With most people having to spend countless hours to commute between home and work, driving has become a source of strain and stress, increasing the possibility of road accidents. Therefore, the research community has been striving to accelerate the shift toward autonomous driving by replacing human drivers with automatic control systems. This shift will improve traffic safety and enhance mobility while boosting human productivity since driving time can be used to do productive tasks instead. The considerable advances in artificial intelligence and information processing technologies further accelerate the shift to autonomous driving. The latter is a complex multidisciplinary process, involving sensing, perception, planning, and control. Control is the final and most important step of the process, it can be divided into longitudinal control, in charge of speed tracking, and lateral control, which handles the steering.

Research works can be divided into two main categories; the first one addresses the longitudinal and lateral controls separately, and the second one couples both tasks together. For instance, Xu \textit{et al.} \cite{xu2018accurate} introduced an optimized controller for speed regulation, where road slope, speed profile, and vehicle dynamics are integrated into the model. Paper \cite{kebbati2021optimized} addressed the longitudinal control by a self-adaptive PID controller, whose optimization and adaptation were based on neural networks and genetic algorithms. An adaptive neural network PID controller was developed by Han \textit{et al.} \cite{s17061244} for the path-tracking task. The authors applied it to a second-order vehicle model, and they used a forgetting factor least square algorithm to estimate model parameters. Guo \textit{et al.} \cite{Guo2019} dealt with path tracking. They developed an MPC controller that takes into account the changing road conditions and small-angle assumptions as a form of measurable disturbance. The authors used the differential evolution algorithm to solve the control problem. Authors of \cite{akbari2020fuzzy} developed a model predictive controller (MPC) for lateral control and optimized its weighting matrices using fuzzy inference systems (FIS). Corno \textit{et al.} \cite{Corno2020} developed an LPV $H_{\infty}$ lateral controller, where they exploited the lateral error and look-ahead distance of the vehicle to ensure better robustness and account for actuator nonlinearities under low speeds. The authors tested their controller in high-speed driving scenarios and evasive maneuvers. Authors of \cite{keb2021} worked on adaptive MPC designed with Laguerre functions for path tracking, they optimized the controller tuning with an improved PSO algorithm. A lookup table approach was used to achieve online controller adaptation. However, the employed method cannot account for all possible cases despite the good results that were achieved. Additionally, the same authors enhanced their approach in \cite{kebb2021} and replaced the lookup table approach with neural networks and adaptive neuro-fuzzy inference systems so that the adaptions generalize beyond the lookup table data. Although significant tracking improvements were achieved, this approach still requires long offline optimizations.

In \cite{kebb2022}, a coordinated lateral and longitudinal control using LPV-MPC for lateral control with PSO-PID for speed regulation was proposed.   In other works, Yao \textit{et al.} \cite{Yao2019} developed an (MPC) path tracking controller that includes longitudinal speed compensation, their approach aims to overcome the assumption of constant longitudinal speed along the control horizon. This technique seeks to minimize the control deviation caused by fast speed and acceleration variations. In \cite{Wang2019b}, Wang \textit{et al.} designed an improved (MPC) control strategy that includes an adaptive fuzzy controller, the latter aims to change the weights of the cost function to tackle the problem of ride discomfort caused by fixed weights in the standard MPC. Li \textit{et al.} proposed in his paper \cite{Li2019a} an (NMPC) for trajectory tracking, their controller was based on nonlinear vehicle dynamics, and Pacejka tire model \cite{Pacejka2008}, and it tracks the yaw angle and lateral position. Most of the above-mentioned studies ensure autonomous driving using coordination between lateral and longitudinal control. However, full autonomy requires handling both lateral and longitudinal controls simultaneously, especially during highly dynamic maneuvers.

Kebbati \textit{et al.} \cite{kebbati2022autonomous} addressed the coupled lateral and longitudinal control, their solution was based on an LPV-MPC approach with genetic algorithm optimization. A Takagi-Sugeno-based MPC (TS-MPC) for autonomous driving was proposed by Alcala \textit{et al.} \cite{alcala2020ts}. This data-driven approach was merely used to learn a Takagi-Sugeno representation of the vehicle dynamics, which was used by the MPC controller with a Moving Horizon Estimator (MHE) to achieve coupled longitudinal and lateral control. Papers \cite{borreli2013,borreli2007,Rick2019} used the nonlinear model predictive control (NMPC) for autonomous driving and parking applications. Kabzan \textit{et al.} \cite{Kabzan2019a} proposed an online learning MPC controller for autonomous racing by learning the model errors online using Gaussian process regression. Similarly, paper \cite{kebbati2023learning} addressed autonomous racing using an online learning NMPC with Gaussian process regression and online moving horizon state estimation. Paper \cite{tuatulea2020design} dealt with autonomous race driving, the authors introduced an (NMPC) for the reduced F1/10 platform. The authors interpolated the circuit boundaries using $3^{rd}$ order polynomials and implemented them as inequality constraints on the lateral position to keep the vehicle inside the track. Kloeser \textit{et al.} \cite{kloeser2020nmpc} proposed an (NMPC) to tackle autonomous racing for a 1:43 scale race car using a singularity-free path parametric model. The authors used partial spatial reformulation of the prediction model to exclude singularities and implemented obstacle avoidance in the optimization problem as a constraint with the objective of maximizing progress on the path. A review of the most widely used control strategies for autonomous driving is provided in \cite{kebbati2022lateral}. 

This paper contributes to the above-mentioned literature by proposing an improved controller for coupled speed and steering control. The main contributions are threefold; First, to ensure real-time application with minimal computing resources and to overcome the heavy computations of NMPC, an adaptive LPV-MPC is developed for autonomous driving. Second, a novel hybrid GA-PSO algorithm is proposed for optimizing the controller's cost function to achieve optimal control actions and automate the controller's tuning process. Third, a Jordan recurrent neural network is designed and trained to learn the tire lateral dynamics by predicting the cornering stiffness coefficients from measurable parameters only, such as velocities and accelerations. 

The article is divided as follows: Section 2 discusses the model of the vehicle's coupled lateral and longitudinal dynamics. Section 3 explains the development of the proposed controller, the adaptation approach using recurrent neural networks, and the optimization of the controller's cost function through the proposed hybrid GA-PSO algorithm. Evaluation results of the learning approach, optimization, and control are presented and analyzed in section 4. Finally, Section 5 provides conclusions and gives perspectives for future work. 

\section{VEHICLE MODELING\label{se:Vehicle modeling}}
\begin{figure}[htb]
\centering
\includegraphics[width=7.5cm,height=4.5cm]{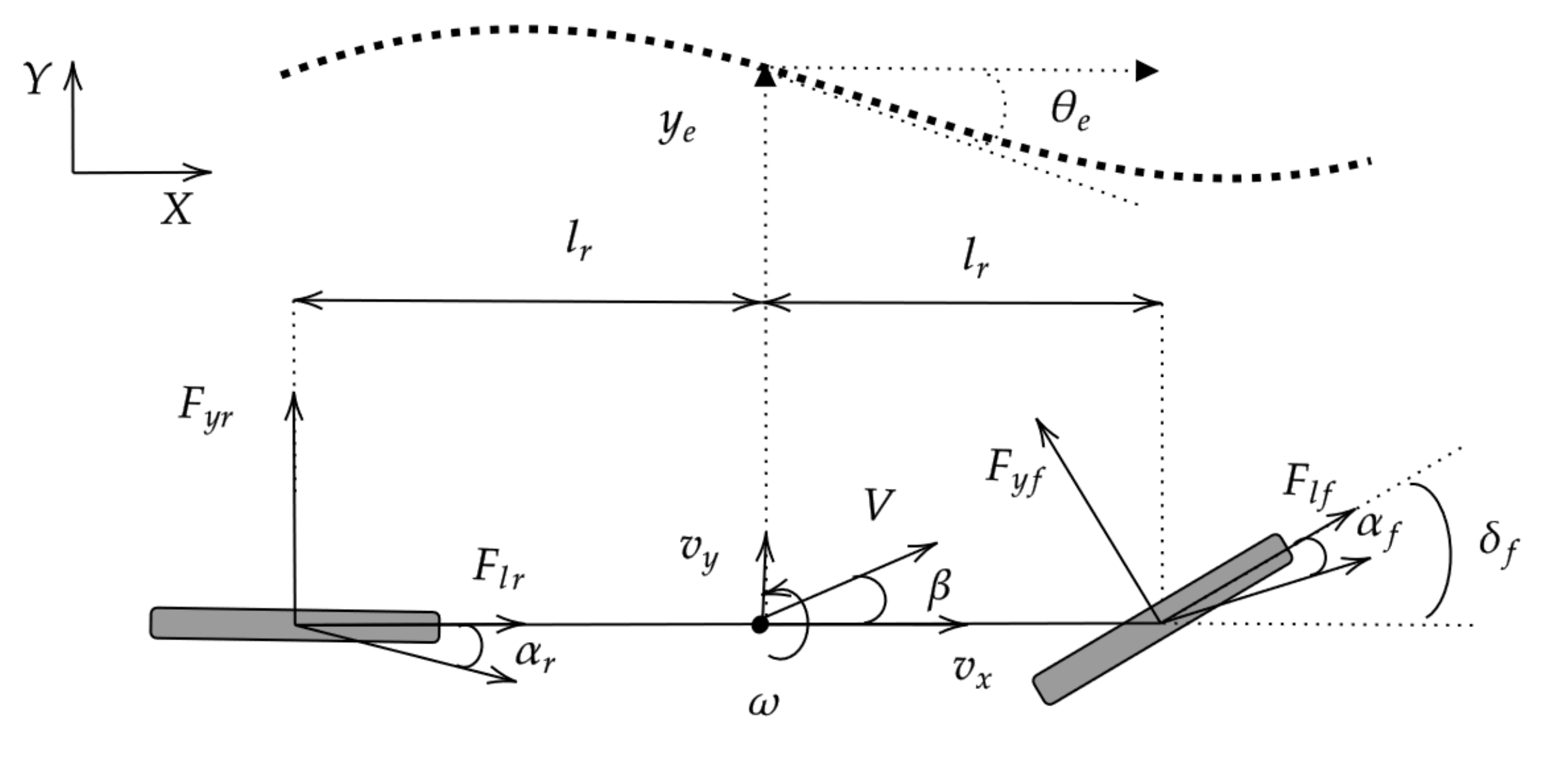}
\caption{Bicycle dynamic model with tracking error.}
\label{fig:1}
\end{figure}

The vehicle is modeled in this article using the common bicycle dynamics model \cite{Alcala2020,schramm2014}, which is considered accurate for control design and easy to implement for real-time control applications as it does not require long computations. Its simplicity is in lumping the front wheels together as well as the rear wheels to form a single-track or bicycle representation as illustrated in Fig. \ref{fig:1}. The tire lateral forces, being a function of the slip angles, govern the vehicle's lateral dynamics. As expressed in equations (\ref{eq1}), the model takes into account the longitudinal, the lateral, and the yaw dynamics and includes the heading and lateral position errors as well:

\begin{equation}
\label{eq1}
\left\{
    \begin{array}{ll}
        \dot{v}_x&= \alpha_x +\omega v_y - \frac{1}{m}(F_{yf}\sin{\delta}+F_d)\\
        \dot{v}_y&= \frac{1}{m}(F_{yf}\cos{\delta}+F_{yr}) -\omega v_x\\
        \dot{\omega}&= \frac{1}{I}(F_{yf}l_f\cos{\delta}-F_{yr}l_r)\\
        \dot{y}_e&= v_x\sin\theta_e + v_y \cos\theta_e\\
        \dot{\theta}_e&= \omega - \frac{v_x\cos\theta_e-v_y\sin\theta_e}{1-y_e k}\\
        F_{yf}&= C_f \alpha_f\\
        F_{yr}&= C_r \alpha_r\\
        F_d&= \mu mg + \frac{1}{2}\rho C_dAv_x^2
    \end{array}
\right.
\end{equation}

The linear longitudinal and lateral velocities and yaw rate in the body frame are represented by $v_x$, $v_y$, and $\omega$, respectively. $F_{y(f,r)}$ express the lateral forces of the front and rear tires, respectively. The total drag force is expressed by $F_d$, where $C_d$, $\rho$, and $A$ represent the drag coefficient, air density, and the vehicle cross-sectional area, respectively. The parameters $\theta_e$ and $y_e$ represent the heading and lateral position errors, where $k$ is the road curvature. The inertia and the mass of the vehicle are represented by $I$ and $m$, and $l_{(f,r)}$ are the distances between the front/rear wheel axles and the vehicle's center of gravity, respectively. The terms $\alpha_x$ and $\delta$ are the acceleration and steering controls, and the parameters $\mu$ and $g$ represent the friction coefficient and the gravity. Finally, the front/rear tire cornering stiffness coefficients are given by $C_{(f,r)}$, and $\alpha_{(f,r)}$ design the slip angles for the front/rear wheels with $\epsilon$ being an additional term to avoid singularities in the model, which are respectively given by:
\begin{equation}
\label{eq4}
\left\{
    \begin{array}{ll}
        \alpha_f&= \delta - \tan^{-1}{( \frac{v_y}{v_x+\epsilon}-\frac{l_f\omega}{v_x+\epsilon})}\\
        \alpha_r&= -\tan^{-1}{ (\frac{v_y}{v_x+\epsilon}+\frac{l_r\omega}{v_x+\epsilon})}
    \end{array}
\right.
\end{equation}
The full model can be considered as a non-linear function that maps the state vector $(x)$ with the input vector $(u)$ and the road curvature $(k)$ as follows:
\begin{equation}
\label{eq5}
\dot{x} = f(x, u, k) 
\end{equation}
where $x = [v_x\ v_y\ \omega\ y_e\ \theta_e]^T$ and $u = [\delta\ \alpha_x]^T$.
\section{CONTROLLER DESIGN}

The control strategy is based on the LPV approach since it allows capturing model nonlinearities, and ensures real-time application with minimal computing resources, therefore,  overcoming the heavy computations of NMPC \cite{Alcala2020}. Furthermore, the LPV approach is adaptive to varying parameters, thereby, increasing the model’s accuracy. The model presented in Section 2 is reformulated in an LPV form and transformed into a state space representation. Thus, the state and control matrices will depend on a scheduling vector of varying parameters. Therefore, the nonlinearities of the model are captured by embedding linear varying parameters into the system matrices, which provides a simple but accurate model for control design. LPV systems are known as a class of linear systems with their parameters being functions of scheduling signals that can be external or internal. The LPV state space formulation of the system is given by (\ref{eq6}) based on the following scheduling vector $\psi=[\delta\ v_x\ v_y\ \theta_e\ y_e\ k ]^T$. 

\begin{equation}
\label{eq6}
\dot{x} = A(\psi)x + B(\psi)u 
\end{equation}
The state matrix $A(\psi)$ and control matrix $B(\psi)$ can then be derived as the following:
\begin{equation}
\label{eq7}
A(\psi) = \left[{\begin{array}{ccccc}A_{11} & A_{12} & A_{13} & 0 & 0\\ 0 & A_{22} & A_{23} & 0 & 0\\0 & A_{32} & A_{33} & 0 & 0\\A_{41} & A_{42} & 0 & 0 & 0 \\ A_{51} & A_{52} & 1 & 0 & 0\end{array}}\right],
\end{equation}
\begin{equation}
\label{eq8}
B(\psi) = \left[\begin{array}{cc} B_{11} & 1 \\ B_{21} & 0 \\ B_{31} & 0\\ 0 & 0 \\ 0 & 0 \end{array}\right],
\end{equation}
where the terms are given by:\\
\noindent
$A_{11} = \frac{-\mu g}{v_x}- \frac{\rho C_d A v_x}{2m},\
A_{12} = \frac{C_f \sin{\delta}}{m v_x},\\
A_{13} = \frac{C_f l_f \sin{\delta}}{m v_x}+ v_y,\
A_{22} = -\frac{C_r + C_f \cos{\delta}}{m v_x},\\
A_{23} = -\frac{C_f l_f \cos{\delta} - C_r l_r}{m v_x}- v_x,\
A_{32} = -\frac{C_f l_f \cos{\delta + C_r l_r}}{I v_x},\\
A_{33} = -\frac{C_f l_f^2 \cos{\delta + C_r l_r^2}}{I v_x},\
A_{41} = \sin{\theta_e},\
A_{42} = \cos{\theta_e},\\
A_{45} = v_x,\
A_{51} = -\frac{k \cos{\theta_e}}{1-y_e k},\
A_{52} = \frac{k \cos{\theta_e}}{1 - y_e k},\\
B_{11} = - \frac{C_f \sin{\delta}}{m},\
B_{21} = -\frac{C_f \cos{\delta}}{m},\
B_{31} = -\frac{C_f l_f \cos{\delta}}{I}.$\\

Generally speaking, the MPC approach exploits the plant model to foresee its behavior over a prediction horizon $N_p$. Based on these predictions, it generates an optimal control sequence by solving a constrained convex optimization problem. The MPC approach is based on the receding horizon principle, where only the first term of the optimal control sequence is used. The LPV model, discretized with $T_s$ sampling time, is used to build the MPC prediction model. This means that the scheduling vector is used to instantiate the LPV model iteratively. The parameters of the scheduling vector can be obtained from sensors, planners, or previous MPC predictions. In this regard, the MPC problem is formulated as the following constrained quadratic optimization:
\begin{equation}
\label{eq9}
\begin{split}
\min_{\Delta U_k} \ J_k & = \sum_{i=0}^{N_p-1} \Big( (r_{k+i} - x_{k+i})^T Q (r_{k+i} - x_{k+i}) + \\
         & \Delta u_{k+i} R \Delta u_{k+i} \Big) + x_{k+N_p}^T Q x_{k+N_p}\\
s.t:\\ 
    & x_{k+i+1} = x_{k+i} + A(\psi_{k+i}) x_{k+i} + B(\psi_{k+i}) u_{k+i}dt \\
    & u_{k+i} = u_{k+i-1} + \Delta u_{k+i}\\
    & \Delta u_{min}  \leq \Delta u_k \leq \Delta u_{max}\\
    &  u_{min}  \leq u_k \leq \ u_{max}\\
    & x_{min}  \leq x_k \leq x_{max}
\end{split}
\end{equation}
The terms $x$, $u$, and $N_p$ define the state vector, the control vector, and the prediction horizon, respectively. The weighting matrices $Q \in \mathbb{R}^{5 \times 5}$ and $R \in \mathbb{R}^{2 \times 2}$ are semi-positive definite, and they penalize the states and the control effort. The longitudinal speed profile is given by the reference vector $r_{k+i}$, and the upper and lower bounds on the control actions, control increments, and states are respectively expressed as $[u_{min}, u_{max}]$, $[\Delta u_{min}, \Delta u_{max}]$ and $[x_{min}, x_{max}]$. The last term of the cost function is added to increase stability, a terminal cost and a terminal set are needed to ensure the asymptotic stability of MPC with a quadratic stage cost. Otherwise, a sufficiently long prediction horizon is required to guarantee asymptotic stability as illustrated by Franz \textit{et al.} \cite{russwurm2021mpc}. In addition, paper \cite{mayne2000constrained} states that adding a terminal cost to the MPC formulation and using a sufficiently long prediction horizon ensures MPC stability without the need for terminal constraints. 

\subsection{Controller Adaptation with Jordan Network} \label{sec:adaptation}

\begin{figure}[htb]
\centering
\includegraphics[width=8.3cm,height=5.25cm]{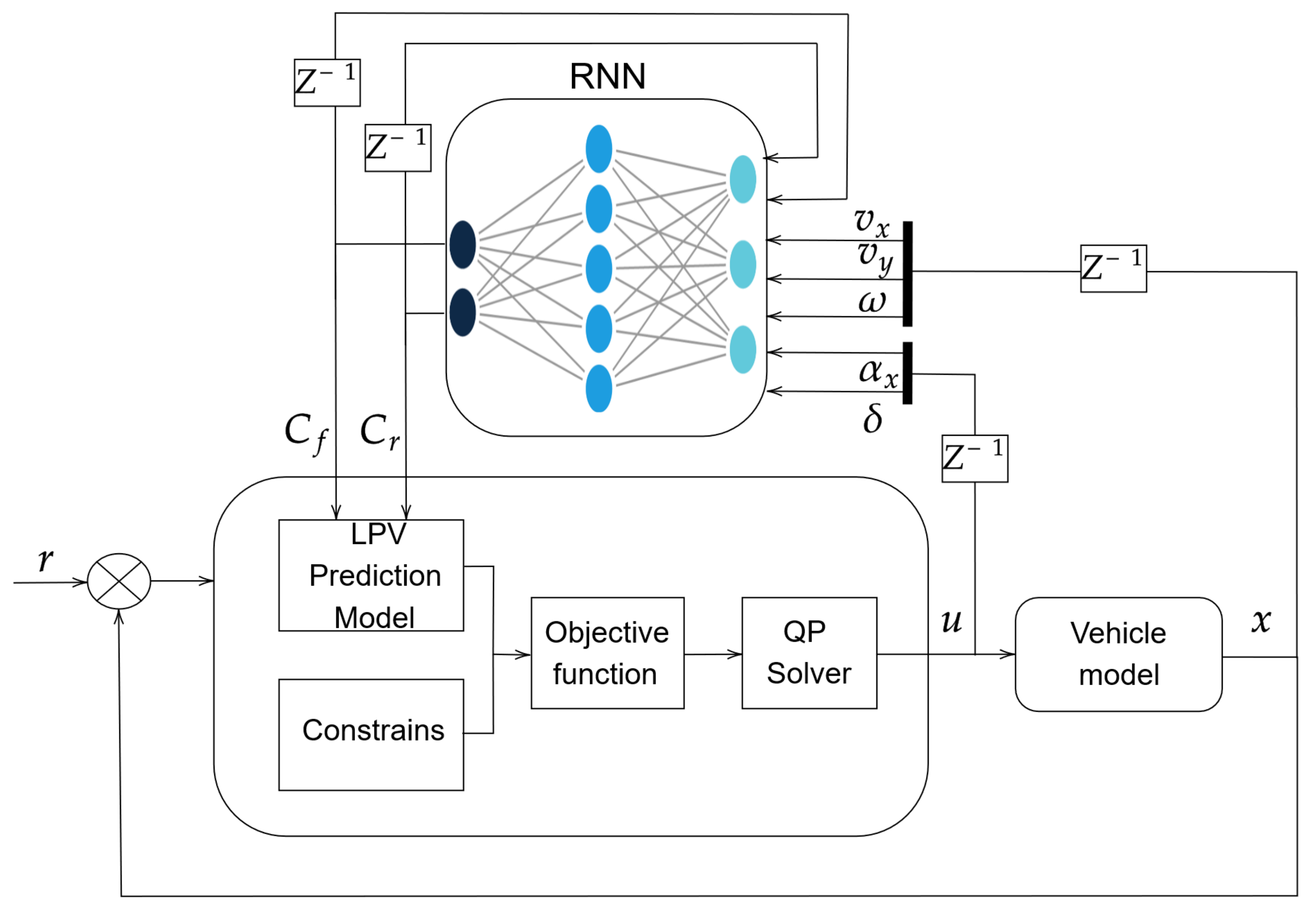}
\caption{Adaptive LPV-MPC approach.}
\label{fig:2}
\end{figure}

\begin{figure}[htb]
\centering
\includegraphics[width=8.3cm,height=4.5cm]{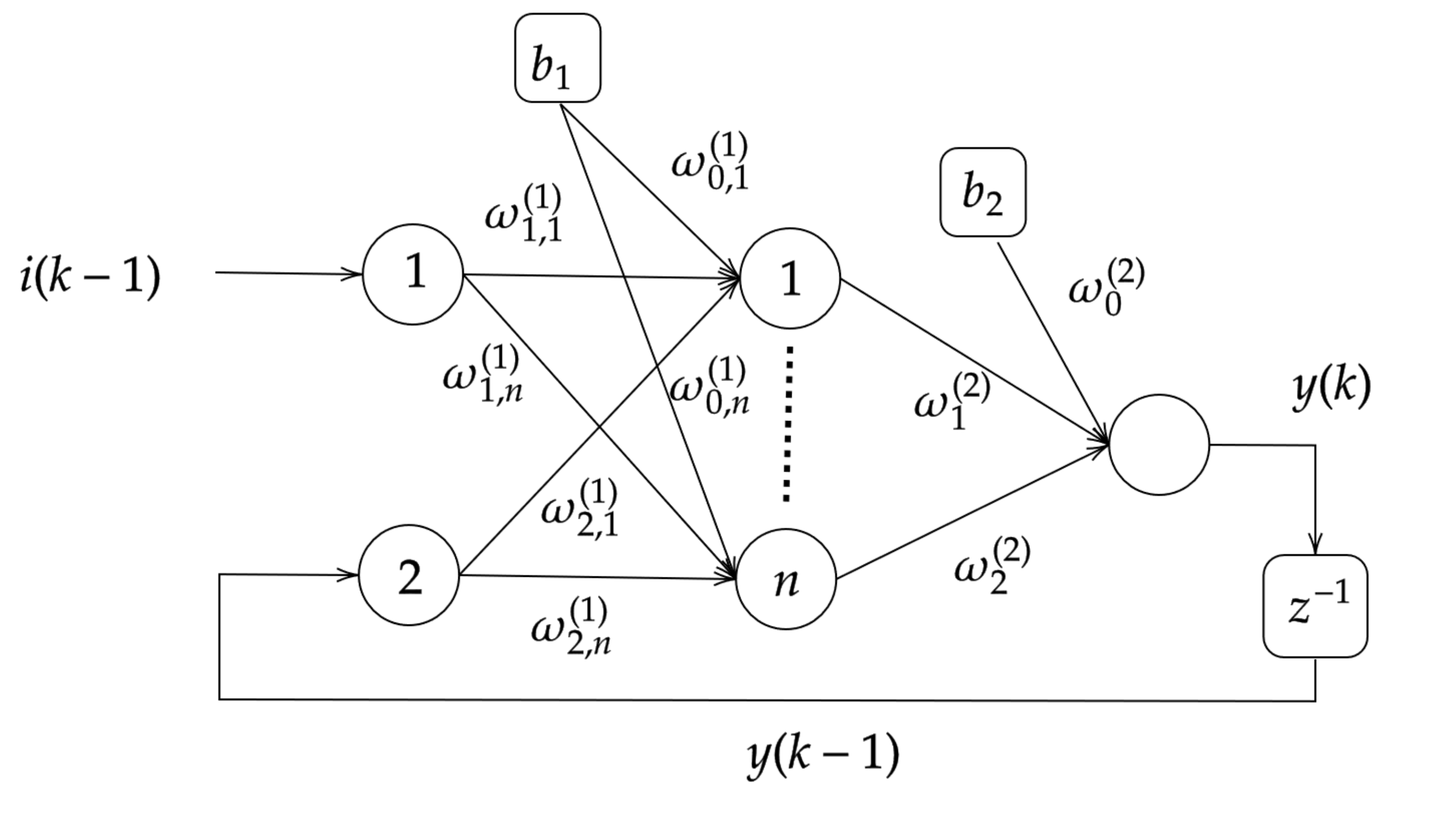}
\caption{Jordan network structure.}
\label{fig:22}
\end{figure}

In most cases, research works are based on linearized tire models, where the tire lateral force depends linearly on the slip angle through a constant known as the cornering stiffness coefficient. However, this is only valid for relatively small slip angles and not during fast and challenging maneuvers, the so-called cornering stiffness coefficient may be time-varying and changes during different types of maneuvers. To deal with this issue, we use machine learning tools to predict the cornering stiffness coefficient online using information from measurable parameters of the vehicle's dynamics. We assume that measurable parameters such as the longitudinal ($v_x$) and lateral ($v_y$) velocities, the steering angle ($\delta$), the acceleration ($\alpha_x$), and the yaw rate ($\omega$)) are enough to capture the tire dynamics. Hence, the prediction model of the LPV-MPC controller is adapted online using this approach (see Fig. \ref{fig:2}), this strategy improves the prediction capability and precision of the LPV-MPC.        
In this article, we propose the use of a modified deep Jordan network to learn tire dynamics. Jordan networks are a simple type of recurrent neural network where the delayed output signal from the output layer is reinjected with the inputs to account for temporal dependencies and improve network predictions, meaning that the previous network predictions become inputs for future predictions. Fig. \ref{fig:22} illustrates the simplest form of a Jordan network where $i$ represents the inputs, $y$ is the output, $\omega$ are the weights, and $b$ is the bias. The network consists of two inputs, one output, and one hidden layer with $n$ hidden neurons. It can be modeled as follows:
\begin{equation}
    y(k) = f(u(k-1), y(k-1))
\end{equation}
Considering $\omega_{j,k}^{(n)}$ as the weights of the $n^{th}$ layer between neurons $j$ and $k$ of the previous and actual layers, respectively, the output signal can be expressed as:
\begin{equation}
    y(k) = \omega_0^{(2)}+ \sum_{i=0}^n \omega_i^{(2)} \zeta(z_i(k))
\end{equation}
The term $\zeta$ is the activation function of the hidden layer, and $z_i(k)$ represents the sum of $i^{th}$ hidden node, and it is given by:
\begin{equation}
   z_i(k)= \omega_{0,i}^{(1)} + \omega_{1,i}^{(1)}u(k-1) + \omega_{2,i}^{(1)}y(k-1)
\end{equation}
  
\subsection{Controller Tuning with Hybrid GA-PSO}

\begin{figure}[htb]
\centering
\includegraphics[width=8.3cm,height=2.55cm]{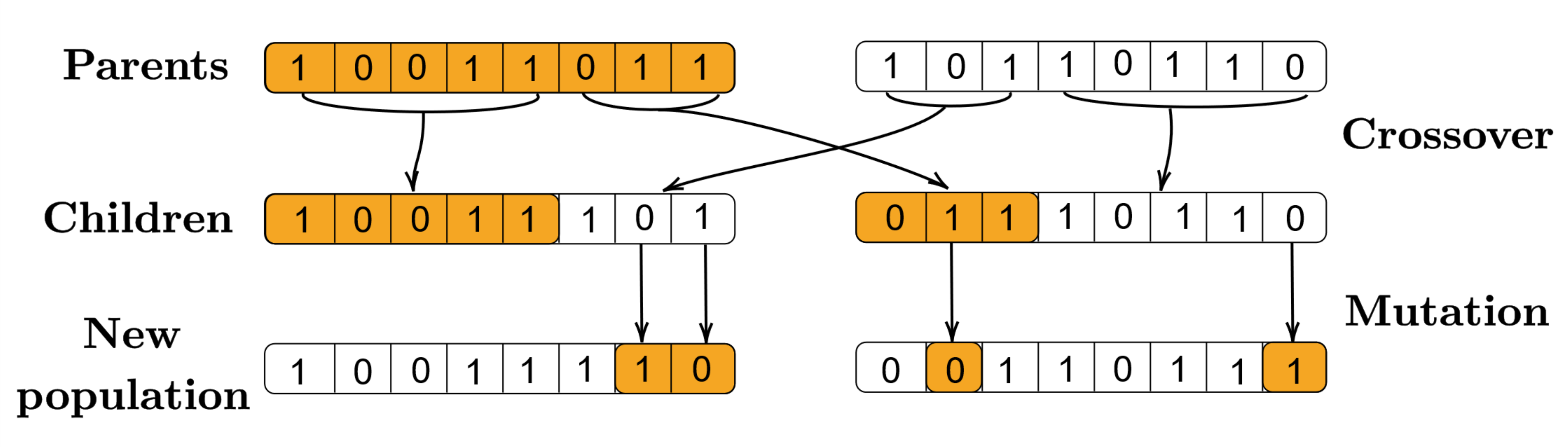}
\caption{Genetic operations}
\label{fig:3}
\end{figure}

Tuning the MPC manually is a difficult task that requires expertise and time, and eventually, it may not result in optimal performance. Thus, to optimize the designed LPV-MPC, we propose a hybrid  Genetic Algorithm (GA) and Particle Swarm Optimizer to tune the weighting matrices of the quadratic cost function by minimizing the MPC tracking root mean squared error (RMSE) as a fitness function. GAs are a global optimization technique based on Darwin’s biological evolution theory. They can find the optimums of discontinuous and nondifferentiable objective functions \cite{song2019}. Generally speaking, the genetic algorithm initializes a population set that contains encoded solutions, which are called chromosomes in the genetic jargon. These possible solutions are improved iteratively and their optimality is assessed by a fitness function. GA evolution operations include the selection, the crossover, and the mutation processes, which control the search capability and the quality of the solutions. They consist of different functions that impact the performance of the algorithm at different degrees \cite{song2019}. For instance, the selection process selects the best chromosomes to be enhanced by the crossover and mutation operations. Alternatively, the crossover seeks to produce high-quality solutions by mixing genetic data, while mutation introduces new genes to complement the crossover as illustrated in Fig. \ref{fig:3}. 
The most common selection operations in the literature, are the roulette wheel (RWS) and tournament selection (TS) \cite{song2019,yadav2017}. For the crossover operation, one finds single/multi-point, uniform, and shuffle crossover. Similarly, mutation operations include inversion and random resetting. Researchers have essentially worked towards improving these operations to further optimize genetic algorithms. In this article, we propose a combination of RWS and TS selection operations to improve the selection of potential genes, a critical GA phase. Briefly speaking, The RWS method provides a higher chance for good genes to be selected, and this improves the exploitation and accelerates the convergence of the algorithm. However, since this approach is mainly based on the fitness value, premature convergence by selecting the same dominant genes is an open issue. On the other hand, the TS approach allows controlling the selection pressure, where smaller tournament sizes ensure more chances for weak genes to be selected unlike RWS. This feature retains the diversity of the search space, which in turn increases the possibility of converging to a global optimum at the expense of slower convergence. In this paper, both methods are used with random percentages at each iteration to increase both convergence speed and optimality by combining the advantages of both methods. Additionally, uniform crossover has been used with mutation based on Gaussian distribution (see Algorithm 1). 
\begin{algorithm}
\caption{Proposed Genetic Algorithm}\label{alg:cap}
\begin{algorithmic}
\Require $Gen_{max},N_p$ \Comment{Generations, Population size}
\State $Pop \gets N_p\ Parents$ \Comment{Random population}
\While{$Generation < Gen_{max}$}
\State $Child \gets empty Pop$ \Comment{Create child population}
    \While{$Child \leq full$}
    \State $RWS \gets \%_r $  \Comment{Generate RWS percentage}
    \State $TS \gets \%_t $
    \If{$\%_r \geq \%_t$}
    \State $Parent1 \gets RWS(Pop)$ \Comment{RWS Selection}
    \State $Parent2 \gets RWS(Pop)$
\Else
    \State $Parent1 \gets TS(Pop)$ \Comment{TS Selection}
    \State $Parent2 \gets TS(Pop)$
\EndIf
\State $Child1,2 \gets UCrossover (Parent1, Parent2)$

\Comment{Perform Uniform Crossover}
\State $Child1,2 \gets GMutation(Child1, Child2)$

\Comment{Perform Mutation}
\State $Fitness \gets Evaluate(Child1, Child2)$

\Comment{Evaluate new offsprings}
\State $Offspring \gets Child1, Child2$

\EndWhile
\State $Pop \gets Offspring$   \Comment{Replace Population}
\EndWhile

\State $Solution \gets Best\ fitness$ \Comment{Save best solution}
\end{algorithmic}
\end{algorithm}

\begin{figure}[htb]
\centering
\includegraphics[width=8.5cm,height=10cm]{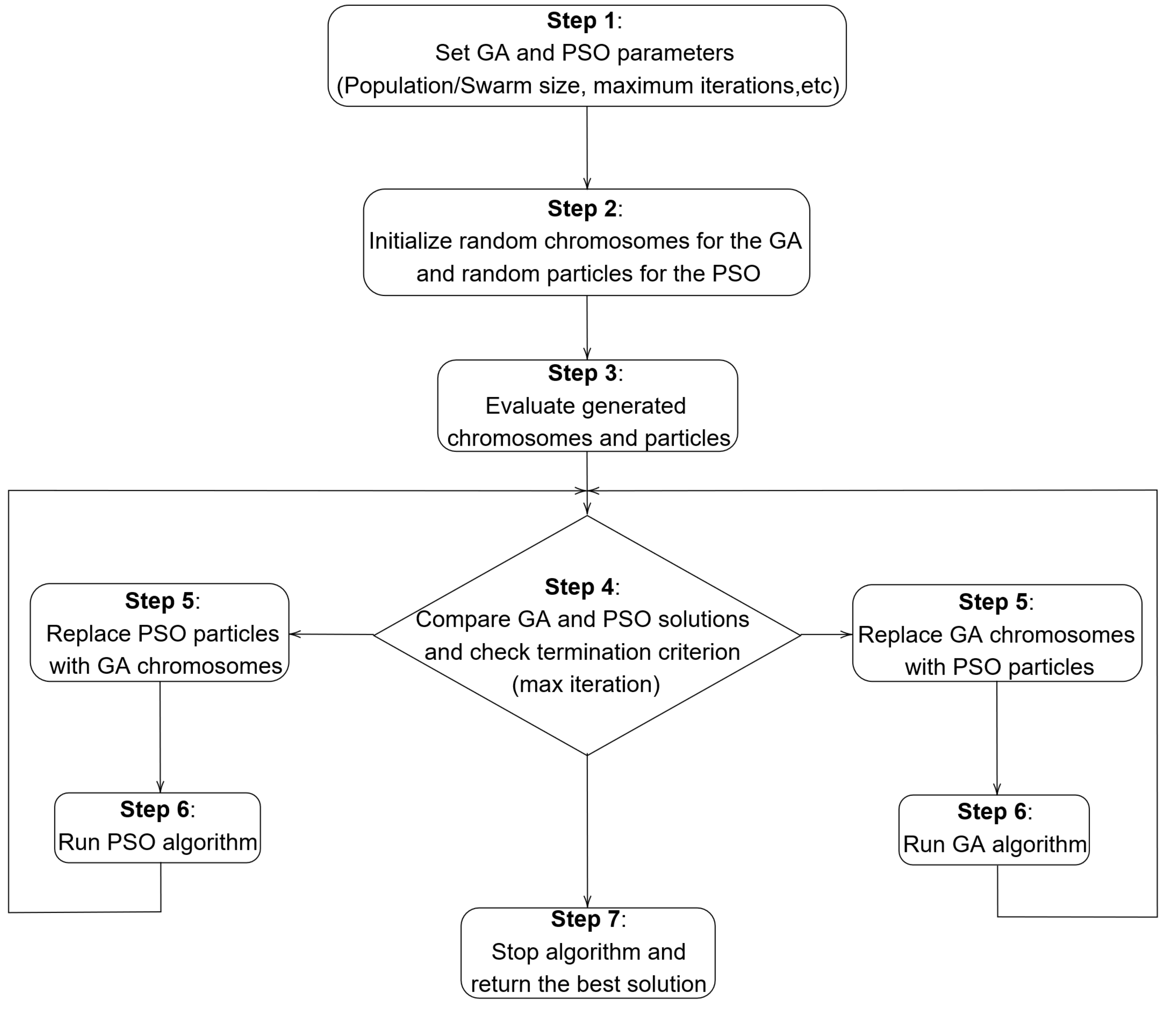}
\caption{Proposed hybrid GA-PSO algorithm.}
\label{fig:12}
\end{figure}

On the other hand, particle swarm optimization is a well-known algorithm for meta-heuristic optimization \cite{song2021improved,xie2019novel}, its classic algorithm is defined as follows:
\begin{equation}
\label{eq18}
\left\{
    \begin{array}{ll}
 v_i(k+1) =& \omega v_i(k) + c_1 r_1 (Pb_i(k)-x_i(k)) \\
 &+ c_2 r_2(Gb(k)-x_i(k))\\
 x_i(k+1)  =& x_i(k) + v_i(k+1)
    \end{array}
\right.
\end{equation}
The terms $v_i$ and $x_i$ define the velocity and the position of particle $(i)$, a particle is a solution to the optimization problem. The terms $\omega$, $c_1$, and $c_2$ are respectively known as inertia weight, cognitive, and social accelerations, while $r_{1,2} \in [0,1]$ are just random constants. Parameters $Pb$ and $Gb$ are the best local and global positions, respectively. In the classic algorithm, $\omega$ and $c_{1,2}$ are constants, but in the improved version of this work, they are dynamic and change according to the following equations \cite{keb2021}: 
\begin{equation}
\label{eq19}
\omega = \omega_{min} + \frac{\exp{(\omega_{max}-\lambda_1(\omega_{max}+\omega_{min})\frac{g}{G})}}{\lambda_2}
\end{equation}
\begin{equation}
\label{eq20}
\left\{
    \begin{array}{ll}
 &c_1(k+1) = c_1(k) + \alpha \\
 &c_2(k+1) = c_2(k) + \beta \\
 &\alpha = -\beta = 0.05\quad \text{for} \quad \frac gG \leq 20\%\\
 &\alpha = -\beta = 0.02\quad \text{for} \quad 20\% \leq \frac gG \leq 35\%\\
  &\alpha = -\beta = -0.035\quad \text{for} \quad 35\% \leq \frac gG \leq 75\%\\
 &\alpha = -\beta = -0.0015\quad \text{for} \quad \frac gG \geq 75\%
    \end{array}
\right.
\end{equation}
The terms $\omega_{min}$ and $\omega_{max}$ are the upper and lower bounds of the inertia weight and $\lambda_{1,2}$ are adjustable parameters to control the decrease from $\omega_{max}$ to $\omega_{min}$. The terms $g$ and $G$ represent the actual and the last generations. The advantage of this improved version over the standard one is that it enhances the overall search capabilities of the PSO algorithm \cite{keb2021}. When $\omega$ decreases exponentially it accelerates the convergence towards the global best solution. Furthermore, increasing cognitive acceleration $c_1$ enhances the exploration phase where particles are pulled towards $Pb$, and increasing $c_2$ enhances the exploitation phase where particles converge towards $Gb$ and vice versa. Compared to GA, PSO algorithms are a bit more intelligent as they incorporate memory by retaining knowledge of good solutions by all the particles as they share information in the swarm. In contrast, a GA would discard all the previous knowledge of the problem once it changes populations.\\
The proposed hybrid algorithm (see Fig. \ref{fig:12}) runs offline and exploits the improved GA and the PSO algorithm iteratively. Briefly speaking, at each iteration of the algorithm, the solutions found by the GA and PSO algorithms are compared, and only the best-found solution is retained. At the next iteration stage, both algorithms will run with the best previously found solution. This strategy allows combining both search efforts of the GA and PSO algorithms towards finding the best solution to the problem. The interest in hybridizing these algorithms lies in the fact that it allows us to overcome the weak searching ability and slow convergence of the GA. In fact, when individuals are not selected in the GA algorithm, their information is completely lost. This is not the case with PSO, since it has memory. On the other hand, PSO algorithms do not have selection operators, which means the algorithm will likely run computations on unfit individuals. Simply put, combining GA with PSO retains the advantages of both algorithms, where the GA excels at reaching the global solution region, and the group search feature of the PSO algorithm boosts the search for exact optimal solutions. 

\section{RESULTS AND DISCUSSION}
The proposed control strategy is tested using a Renault Zoe vehicle, whose behavior is simulated in \texttt{Matlab} using a high fidelity nonlinear dynamic model \cite{Alcala2020} with the Pacejka formula for the lateral tire forces \cite{Pacejka2008}. Table \ref{tab:2} properly lists all model parameters.\\

\subsection{Learning and Optimization Results}
\begin{figure}[htb]
\centering
\includegraphics[width=8.35cm,height=5.5cm]{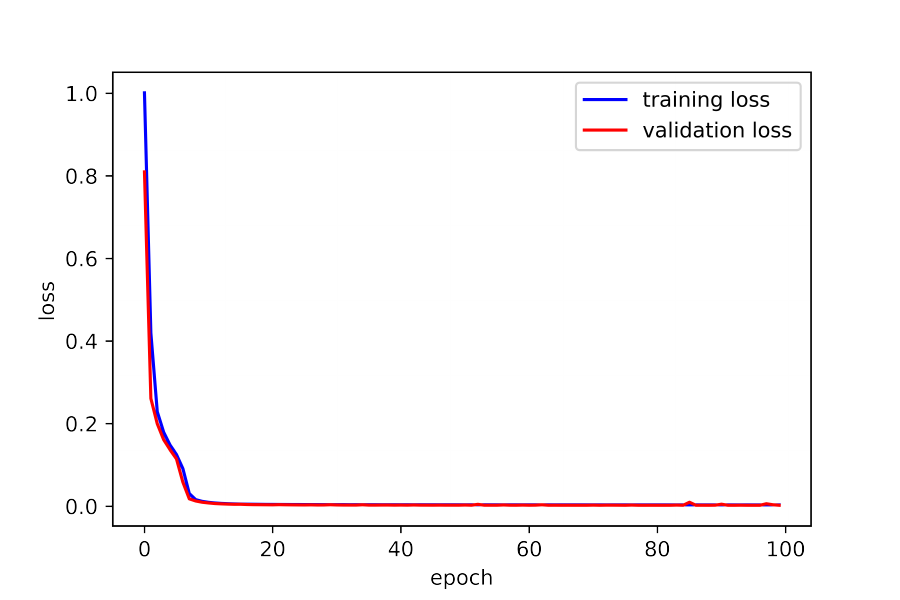}
\caption{Learning curve.}
\label{fig:4}
\end{figure}

\begin{figure}[htb]
\centering
\includegraphics[width=8cm,height=10cm]{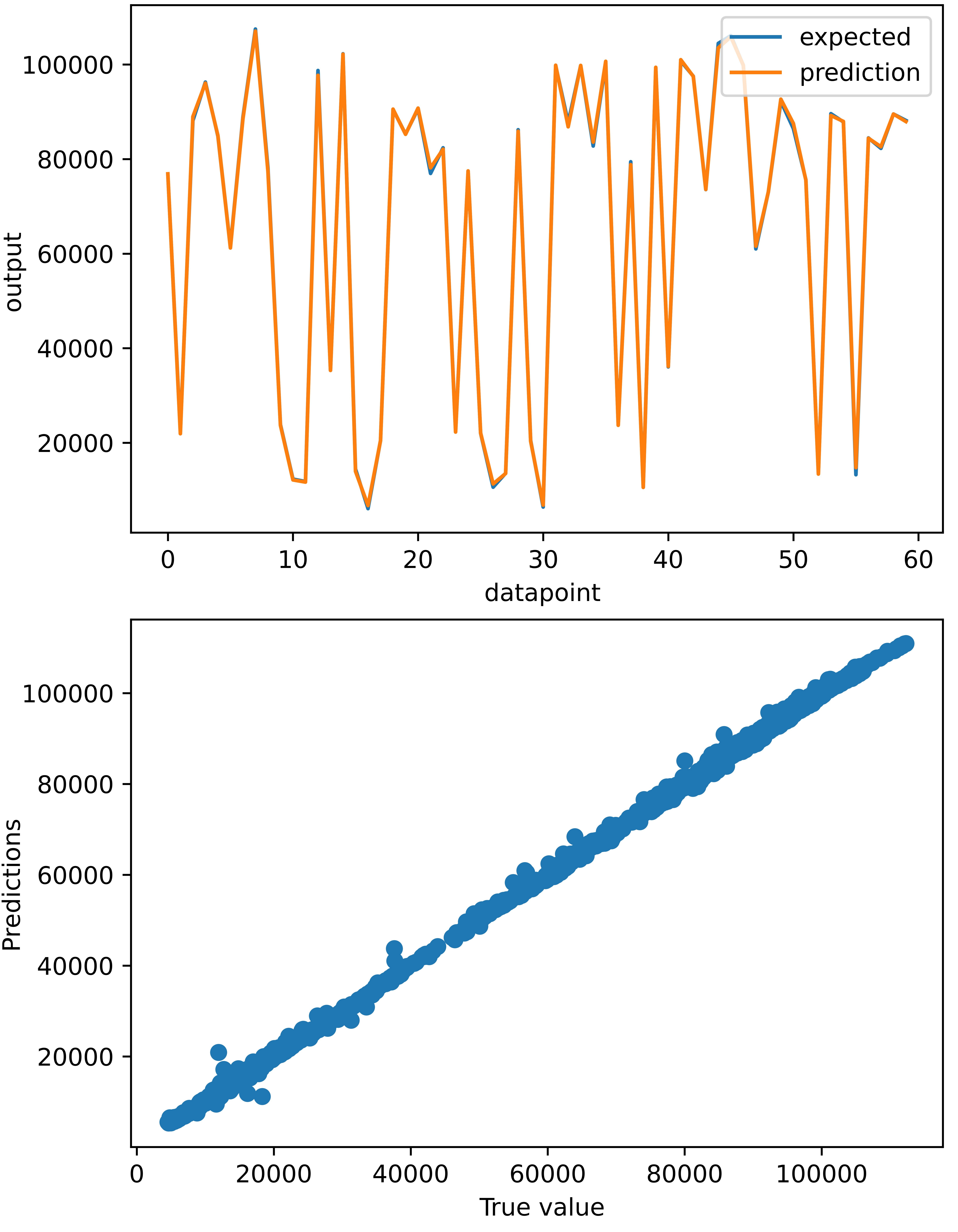}
\caption{Prediction model performance.}
\label{fig:5}
\end{figure}

Since the cornering stiffness coefficients are learned from data, multiple driving scenarios are performed in Carsim \cite{benekohal1988carsim} to collect data for training the neural network. The Jordan recurrent neural network consists of an input layer with $7$ neurons, and two hidden layers with $8$ and $5$ hidden neurons, respectively. In addition to one output layer with two neurons corresponding to the cornering stiffness coefficients for front and rear wheels. All the hidden layers are activated with the sigmoid function, while the identity is used for the output layer. The model is developed in Keras-Tensorflow \cite{ketkar2017introduction} and trained for $100$ epochs with a batch size of $16$. The data set contains around $10000$ data points, which were split into $75\%$ for training and $25\%$ for validation. The adaptive movement estimation algorithm was used as the optimizer with a learning rate of $5e^{-4}$. Fig. \ref{fig:4} shows the training curve of the neural network, which proves the ability of the model to learn the data without over-fitting. The resulting validation loss value is as low as $0.002$, while the training loss reached $0.005$. The obtained $R^2$ scores on the training and the test data sets are $99.8\%$ and $99.7\%$, respectively. Fig. \ref{fig:5} shows a comparison between the predicted values and the expected ones over a few data points of the test data set, which illustrates the accuracy of the model.


\begin{figure}[htb]
\centering
\includegraphics[width=8cm,height=6cm]{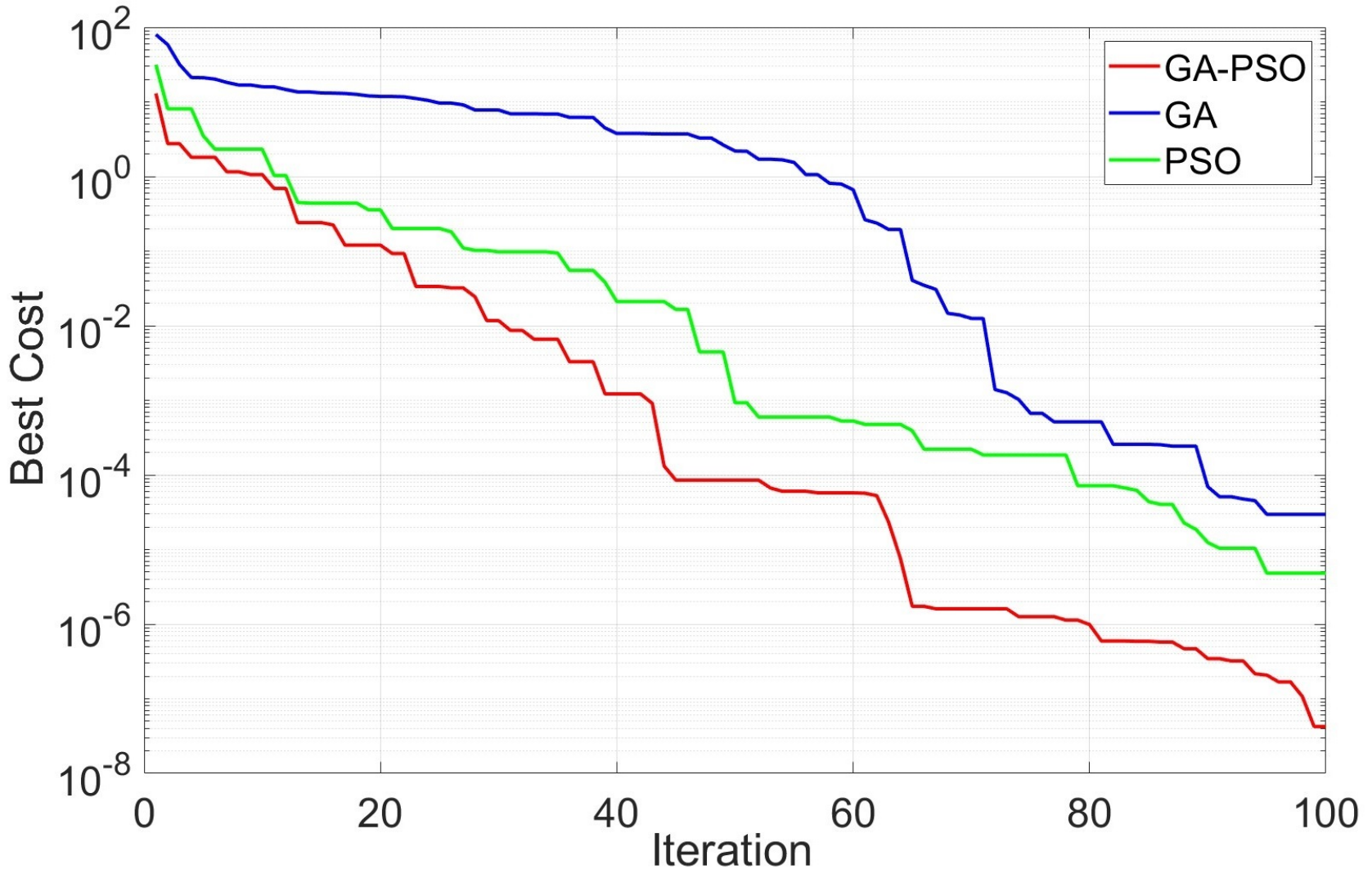}
\caption{Performance of the improved hybrid GA-PSO.}
\label{fig:6}
\end{figure}

The different parameters of the proposed hybrid GA-PSO algorithm are tuned intuitively and iteratively until the desired performance is achieved, the algorithm is evaluated on a 5D sphere function \big($f(x) = \sum_{i=1}^5 x_i^2$\big) as a benchmark test \cite{kaya2011}. The resulting performance of the proposed GA-PSO algorithm over $100$ iterations is compared to improved GA and improved PSO, respectively. Fig. \ref{fig:6} shows that the hybrid GA-PSO is indeed faster and able to further optimize the solutions. It managed to reach a minimum cost value of $4.24e^{-8}$ compared to $4.82e^{-6}$ and $2.95e^{-5}$ for the improved PSO and GA algorithms, respectively. Tables \ref{tab:1} and \ref{tab:11}  list the parameters used in the GA-PSO algorithm for the optimization of the LPV-MPC controller. The fitness function, in this case, was selected as RMSE for longitudinal velocity, lateral position, and heading tracking. The GA-PSO optimization managed to achieve a minimum RMSE score of $0.0069$, $0.0191$, and $0.0212$ for the position, heading, and velocity tracking, respectively. The optimized weighting matrices are as follows:\\
$Q = \text{diag}(50,\ 1e^{-5},\ 0.01,\ 47.43,\ 1e^{-3})^T$,\\
$R = \text{diag}(0.003,\ 1e^{-4})^T$.

\begin{table}[htb]
\caption{GA parameters} 
\label{tab:1}
\centering
\begin{tabular}{c c c} 
\hline
Parameter & Name & Value\\[0.8ex] 
\hline
$Gen$ & Generation  & $15$  \\[0.8ex] 
$N_P$ & Size of population & $25$  \\[0.8ex] 
$O_p$ & Percentage of offsprings & $0.8$  \\ [0.8ex] 
$\beta$ & Selection pressure &  $0.75$ \\ [0.8ex] 
$\mu_r$ & Mutation rate &  $0.3$ \\ [0.8ex] 
$\sigma$ & Mutation variance &$0.15$ \\ [0.8ex] 
\hline 
\end{tabular}
\end{table}

\begin{table}[htb]
\caption{PSO hyper-parameters.} 
\label{tab:11}
\centering
\begin{tabular}{c c c} 
\hline
Parameter & Interpretation & Value\\[0.8ex] 
\hline
$N$ & Generation & $15$ \\[0.8ex] 
$N_{Pop}$ &  Swarm particles & $25$ \\[0.8ex] 
$\omega_{max}$ & Maximum inertia weight & $0.99$ \\ [0.8ex] 
$\omega_{min}$ & Minimum inertia weight & $0.1$ \\ [0.8ex] 
$c_{1i}$ & Initial cognitive acceleration& $2$ \\ [0.8ex] 
$c_{2i}$ & Initial social acceleration& $2$ \\ [0.8ex] 
$\lambda_1$ & Constant & $30$ \\ [0.8ex]
$\lambda_2$ & Constant & $3$ \\ [0.8ex]
\hline 
\end{tabular}
\end{table}
\subsection{Control Results}
\begin{figure}[htb]
\centering
\includegraphics[width=8cm,height=6.5cm]{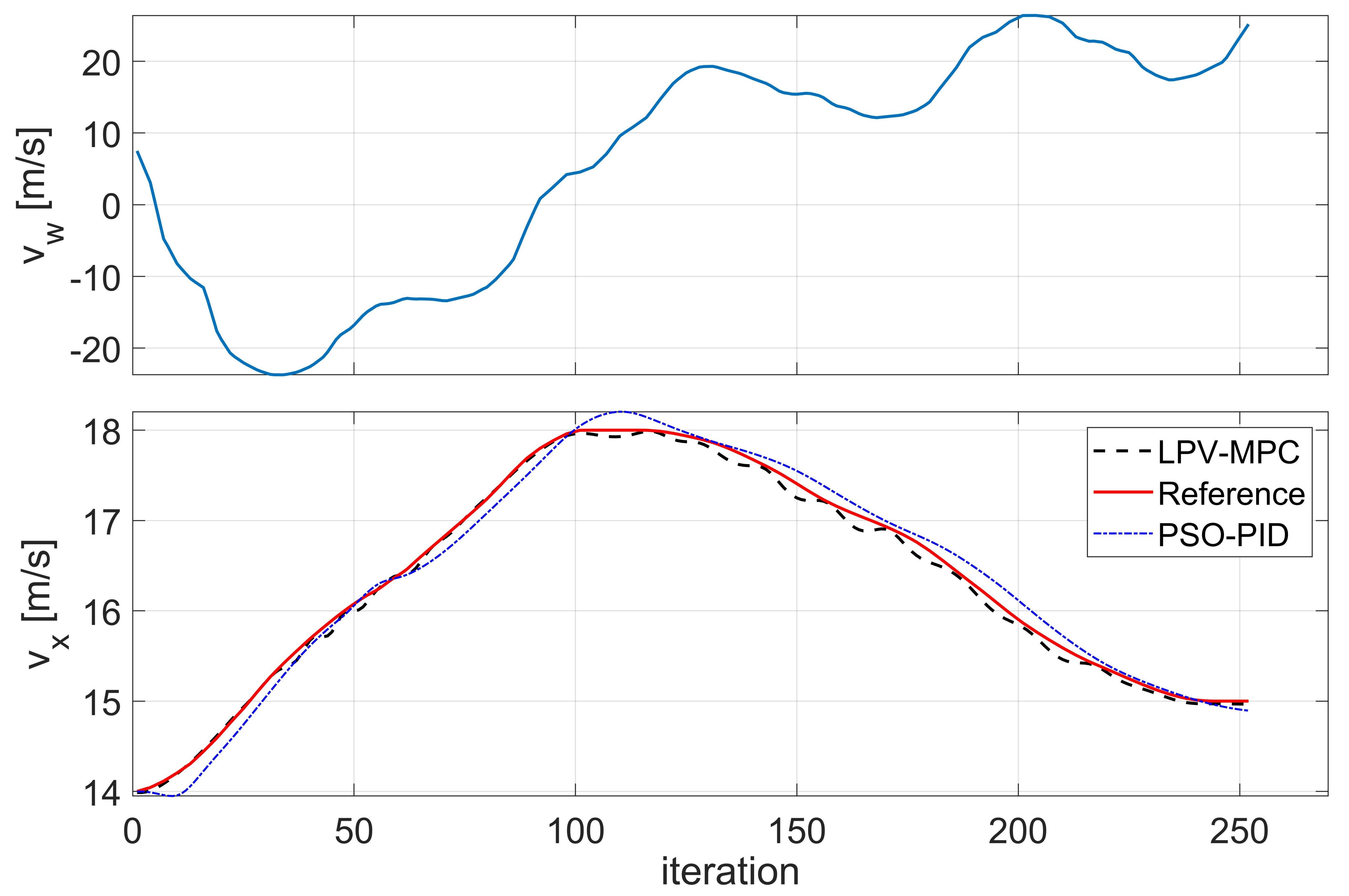}
\caption{Wind velocity and speed tracking for double lane change.}
\label{fig:77}
\end{figure}

\begin{figure}[htb]
\centering
\includegraphics[width=8cm,height=6.5cm]{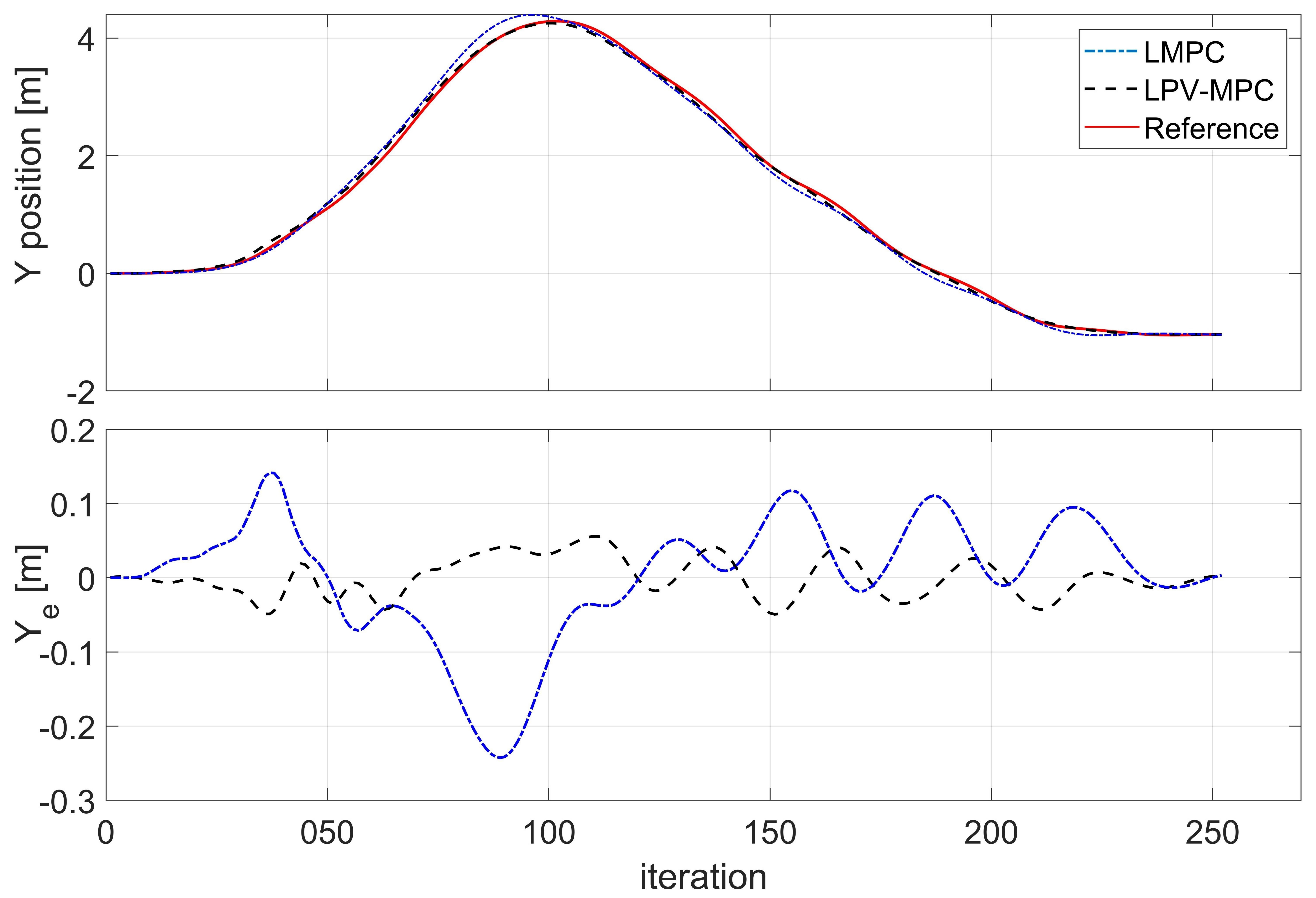}
\caption{Trajectory tracking for double lane change.}
\label{fig:777}
\end{figure}

\begin{figure}[htb]
\centering
\includegraphics[width=8cm,height=6.5cm]{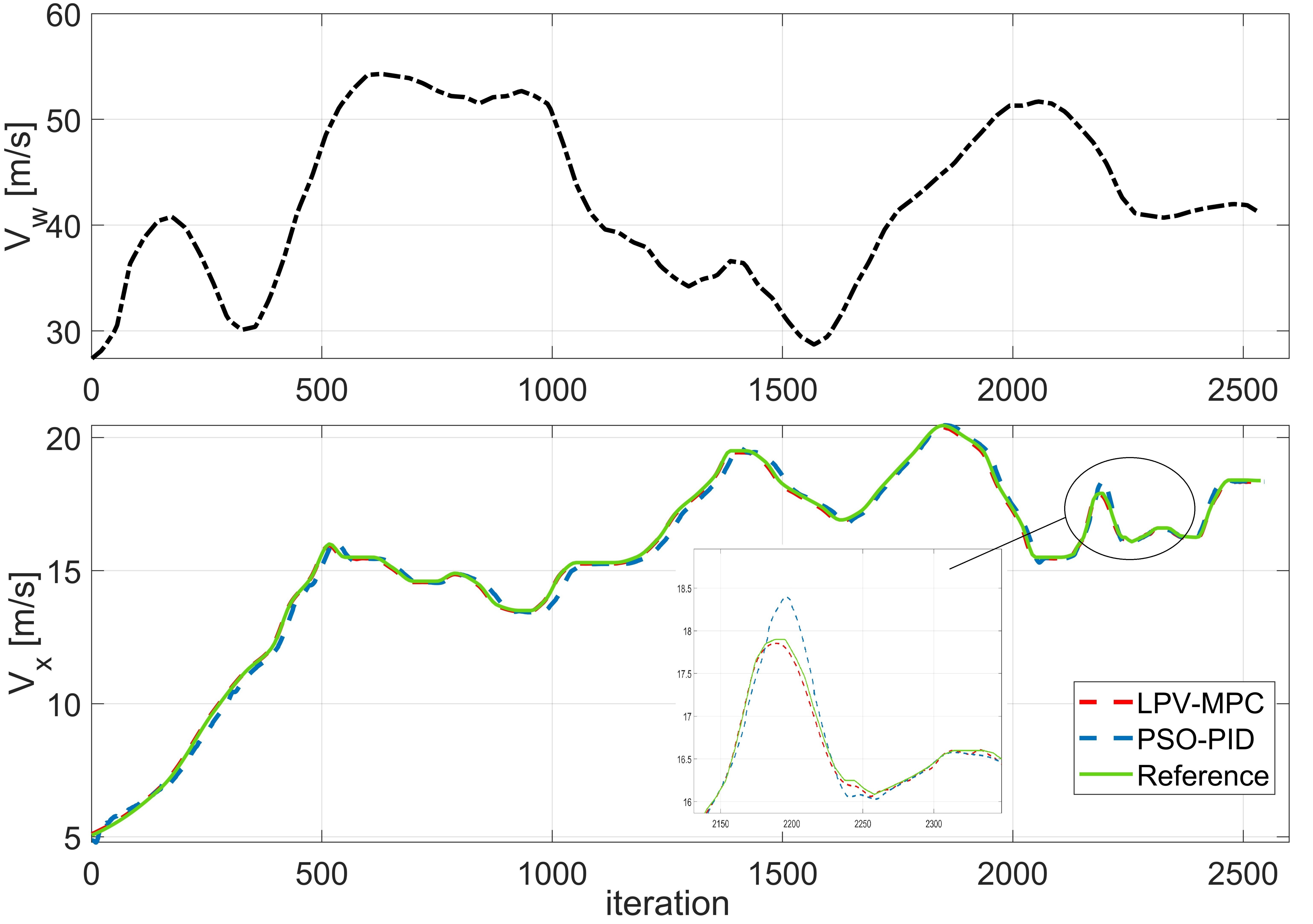}
\caption{Wind velocity and speed tracking performance.}
\label{fig:7}
\end{figure}

\begin{figure}[htb]
\centering
\includegraphics[width=8cm,height=6.5cm]{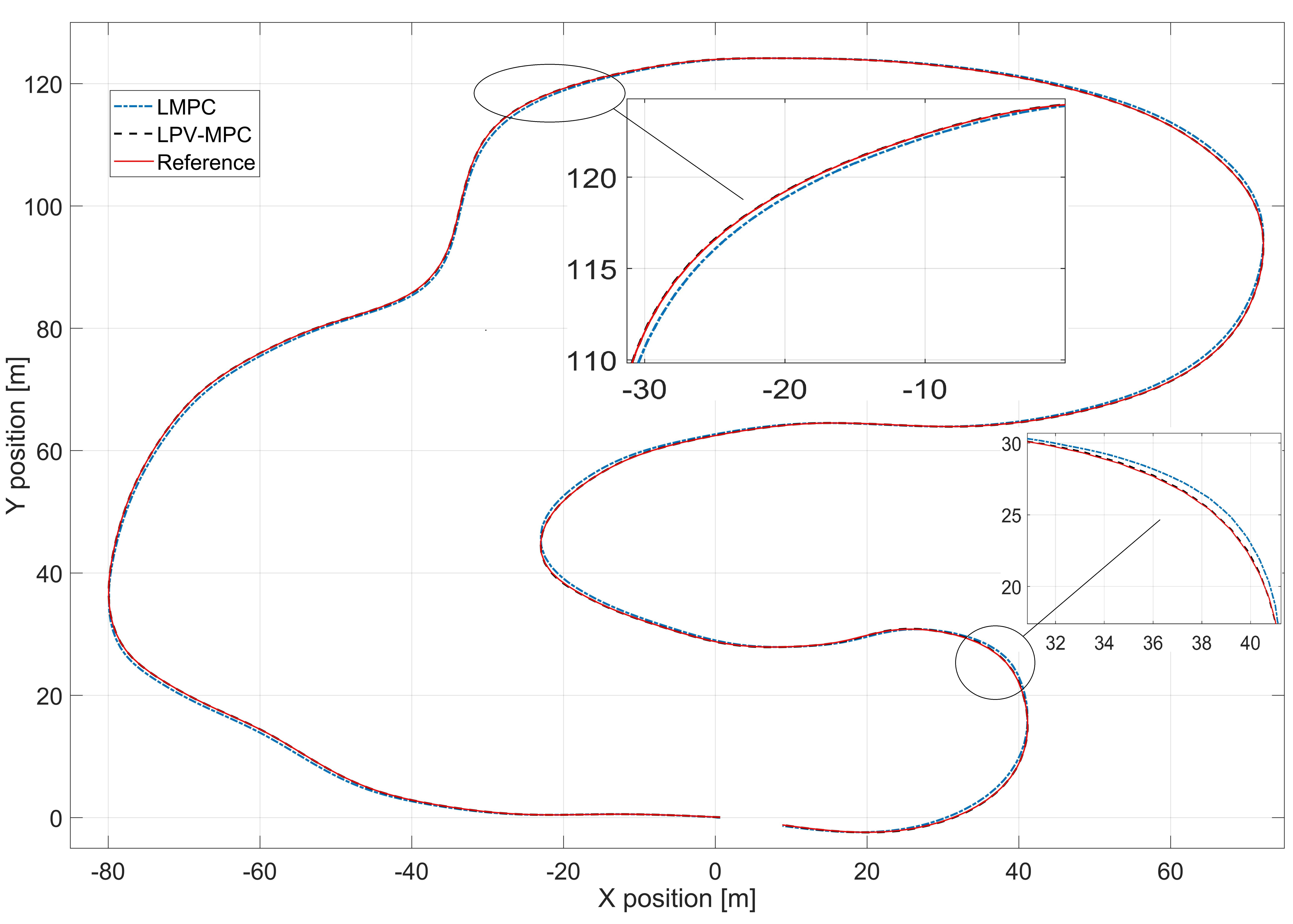}
\caption{Trajectory tracking.}
\label{fig:8}
\end{figure}

The proposed controller is coded in \texttt{Yalmip} platform and solved using \texttt{Gurobi} solver. The algorithm runs at 95Hz on a Ryzen$7$ laptop with $32gb$ of RAM. The LPV-MPC is implemented and evaluated in \texttt{Matlab} simulations using the high fidelity nonlinear dynamic model \cite{Alcala2020} with the Pacejka formula for the lateral tire forces \cite{Pacejka2008}. Table \ref{tab:2} presents the MPC parameters. The evaluation is performed for a double lane change trajectory and speed profile (see Fig. \ref{fig:77},\ref{fig:777}). Overall, the proposed controller performs the double lane change maneuver very well with  minimal tracking errors. Similarly, the controller is further tested on a more challenging general trajectory and speed profile under wind disturbances varying between $20$ and $50\ m/s$, (see Fig. \ref{fig:7},\ref{fig:8}). Furthermore, it is compared to another MPC based on the linear bicycle model as introduced in \cite{kebb2022}, which is used in coordination with an optimized PSO-PID to address the combined longitudinal and lateral dynamics. We denote this controller LMPC for linear MPC. Such a comparison shows that the proposed controller handles both lateral and longitudinal dynamics and outperforms the decoupled control strategy of \cite{kebb2022}, which dedicates two optimized controllers for the same task. The LMPC was tested on the same trajectory with the same speed profile and wind disturbance. Fig. \ref{fig:7} shows the velocity profile variying between $5$ and $25\ m/s$. It can be seen that the LPV-MPC is slightly more accurate in speed tracking. The MSE was evaluated at $0.02$ compared to $0.19$ for the PSO-PID, whose parameters were already optimized. In addition, PSO-PID was found to be more aggressive as it cannot handle constraints.

\begin{figure}[H]
\centering
\includegraphics[width=8cm,height=6cm]{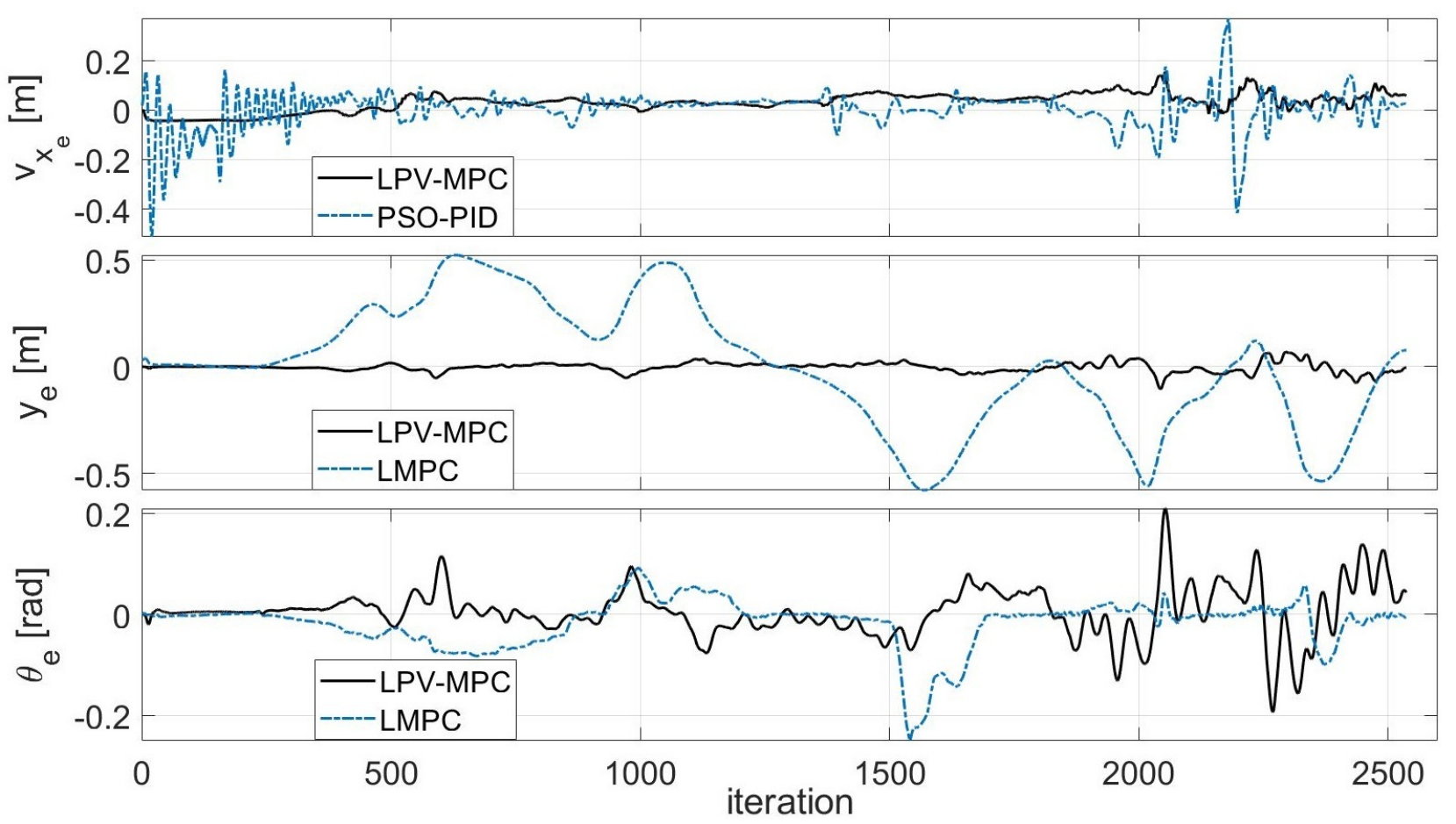}
\caption{Tracking performance.}
\label{fig:10}
\end{figure}

\begin{figure}[H]
\centering
\includegraphics[width=8cm,height=6cm]{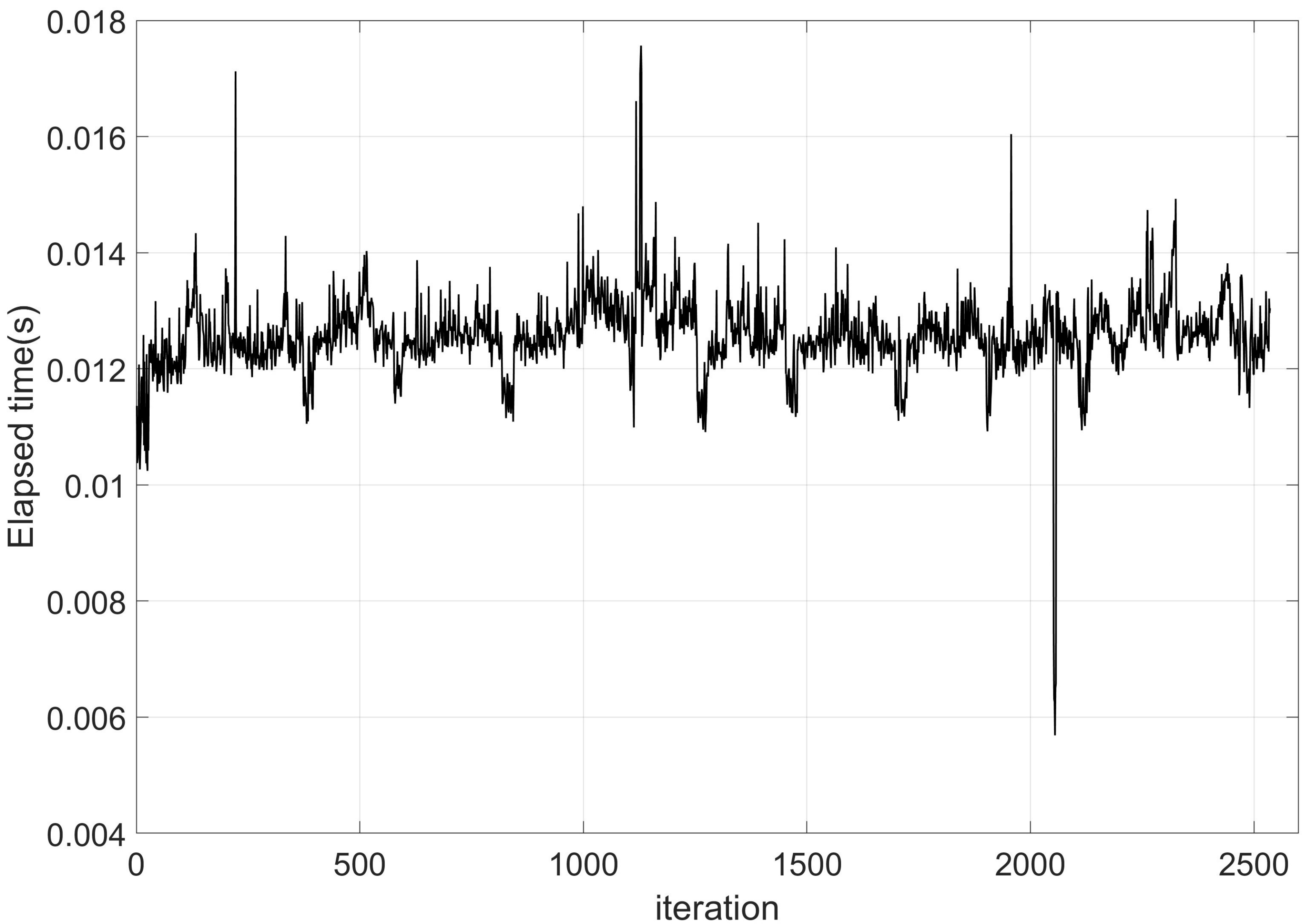}
\caption{MPC computation time.}
\label{fig:11}
\end{figure}

\begin{table}[htb]
\caption{MPC and model parameters} 
\label{tab:2}
\centering
\begin{tabular}{c c c c} 
\hline
Parameter & Value & Parameter & Value\\[0.8ex] 
\hline
$m$ & $1575\ (kg)$ & $C_d$ & $0.29$  \\[0.8ex] 
$I_z$ & $2875\ (kg.m^2)$ & $A$ & $1.6\ (m^2)$  \\[0.8ex] 
$l_f$ & $1.2\ (m)$ & ${y_e}_{max/min}$ & $0.3\ (m)$ \\ [0.8ex] 
$l_r$ & $1.6\ (m)$ & $u_{max/min}$ & $\pm \frac{\pi}{6}\ (rad)$ \\ [0.8ex] 
$\rho$ & $1.225\ (kg m^3)$ & $\Delta u_{max/min}$ & $\pm \frac{\pi}{12}\ (rad)$ \\ [0.8ex] 
$\mu$ & $0.82$ & $N_p$ & $10$  \\ [0.8ex] 
$g$ & $9.81\ (m/s^2)$  & $T_s$ & $0.033\ s$  \\ [0.8ex] 
\hline 
\end{tabular}
\end{table}
The results in Fig. \ref{fig:8} show that the proposed LPV-MPC performed much better than LMPC in terms of tracking accuracy. The obtained MSE score was as low as $0.007$ compared to $0.32$ for LMPC, this means the LPV-MPC tracking is almost ideal compared to its rival. The zoomed regions of the figure further illustrate the big difference in tracking accuracy, this is partly due to the different models used to develop the MPC controller. 
\begin{figure}[htb]
\centering
\includegraphics[width=7.75cm,height=6.5cm]{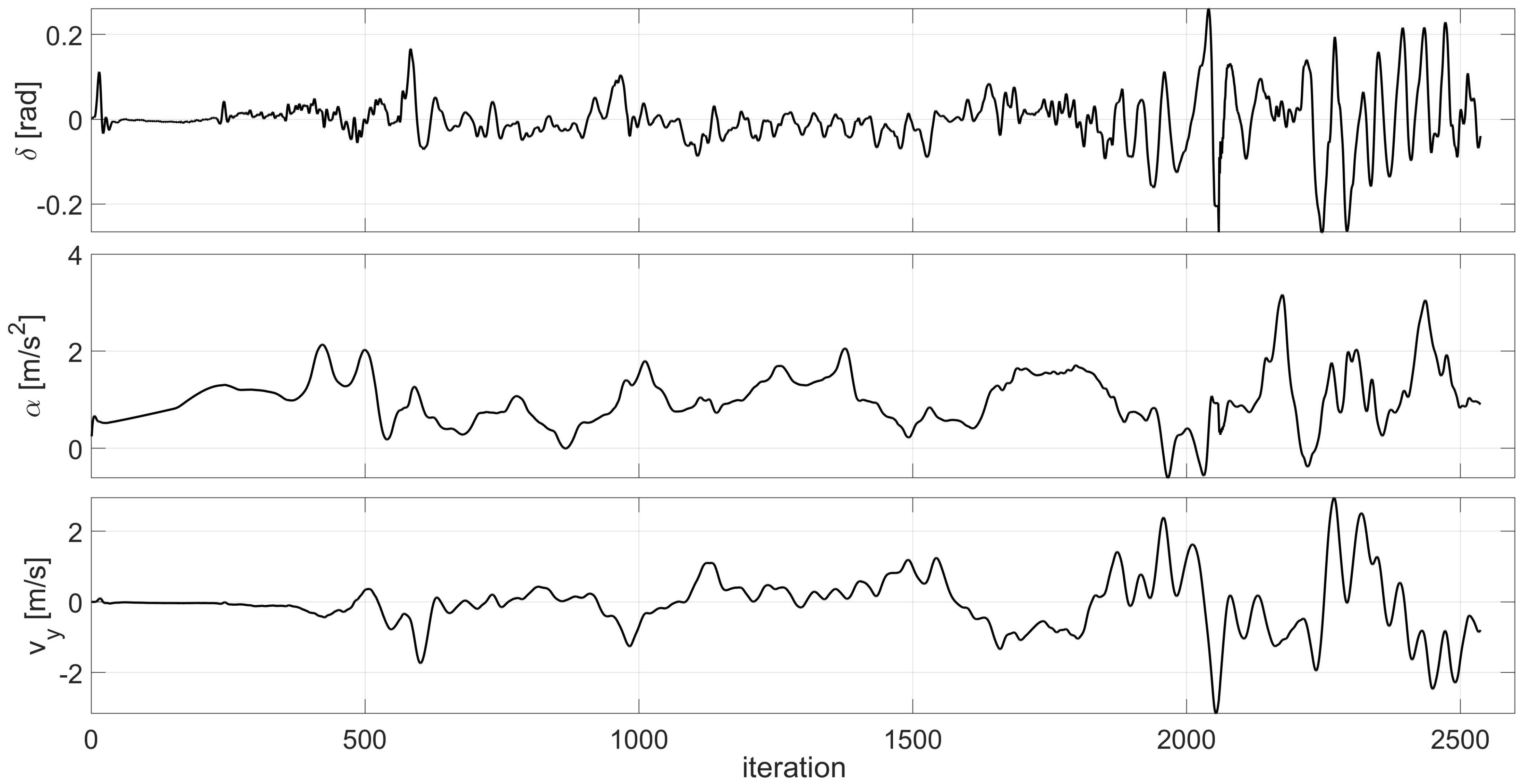}
\caption{Steering, acceleration, and lateral velocity signals.}
\label{fig:9}
\end{figure}
\begin{figure}[H]
\centering
\includegraphics[width=8cm,height=4.25cm]{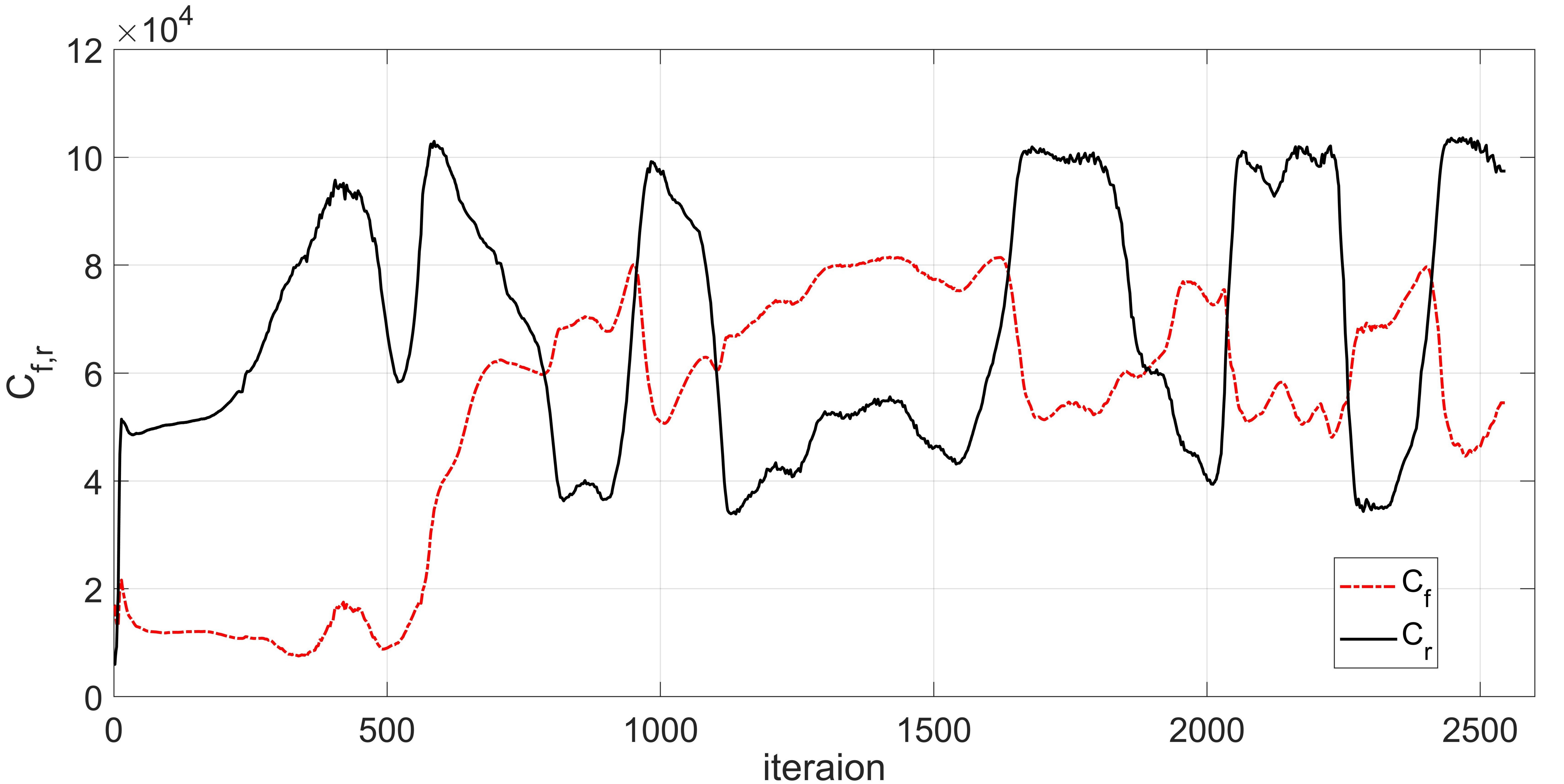}
\caption{Cornering stiffness coefficients.}
\label{fig:100}
\end{figure}
The corresponding steering and acceleration controls and the lateral velocity of the vehicle are shown in Fig. \ref{fig:9}. The predicted cornering stiffness coefficients for the tested trajectory are reported in Fig. \ref{fig:100}, and Fig. \ref{fig:10} shows the corresponding tracking errors for longitudinal velocity, heading, and lateral position, respectively. As seen in the figure, the velocity tracking error does not exceed $0.095\ m/s$, and the maximum heading and position tracking errors are kept below $4.5^{\circ}$ and $2.3\ cm$, respectively. Moreover, the execution time of the LPV-MPC is very suitable for real-time applications as observed in Fig. \ref{fig:11} with a mean computation time of $(0.011\ s)$.

\section{CONCLUSIONS}

This article addressed the coupled control task in autonomous driving with an LPV-MPC controller. The developed controller is capable of simultaneously controlling the lateral and longitudinal dynamics. A machine learning approach has been introduced to predict the model's tire cornering stiffness coefficients online, using only measurable parameters. This approach adapts the LPV-MPC prediction model for more accurate predictions. For tuning and optimizing the proposed controller, an improved hybrid GA-PSO algorithm has been proposed. The developed controller has been evaluated on a challenging track and compared to another variant of LPV-MPC. The obtained results showed superior performance of the proposed controller, which ensures high speed and trajectory tracking accuracy. Future works shall target online learning for more advanced autonomous driving applications.

\section{ACKNOWLEDGEMENTS}
The authors would like to thank the University of Paris-Saclay for the financial support provided to conduct this research.

\section{AUTHORS' CONTRIBUTIONS}
\textbf{Yassine Kebbati}: Conceptualization, Methodology, Writing original draft, Reviewing and Editing; \textbf{Naima Ait-Oufroukh}: Review and Validation; \textbf{Dalil Ichalal}: Review and Supervision; \textbf{Vincent Vigneron}: Supervision, Review and Validation.

\section{CONFLICTS of INTEREST}
The authors declare having no conflict of interest for the publication of this article.

\bibliographystyle{ieeetr}
\bibliography{PhD}

\end{document}